\documentclass{article}

\usepackage{amssymb}
\usepackage[cp1251]{inputenc}
\usepackage[english]{babel}
\usepackage{amsmath}
\usepackage{color}
\usepackage{graphics}
\usepackage{epsfig}
\usepackage[all,knot]{xy}
\xyoption{arc}

\def\Inn{{\bf Inn}}
\def\rg{{\frak g}}
\def\Aut{{\hbox{\bf Aut}\;}}
\def\Out{{\hbox{\bf Out}\;}}
\def\mapr#1{\smash{\mathop{\buildrel{#1}\over\longrightarrow}}}
\def\cD{{\cal D}}
\def\Der{{\hbox{\bf Der}\;}}
\def\cO{{\cal O}}
\newtheorem{theorem}{Theorem}
\def\g{{\bf g}}
\def\Int{{\hbox{\bf Int}\;}}
\def\GL{{\;\hbox{\bf GL}\;}}
\def\Hom{{\hbox{\bf Hom}\;}}
\def\ad{{\hbox{\bf ad}}}

\author{A.S.Mishchenko \\(Moscow Lomonosov State University), \\Leanh Nguyen \\(National University of Civil Engineering, Vietnam)}
\title{Some results on the Mackenzie obstruction for transitive Lie algebroids}

\begin{document}
\maketitle

\begin{abstract}
The preprint is prepared as description of results that were obtained during
joint scientific project No: 71NC /2015/VNCCCT on the VIASM (Vietnam Institute for Advanced Study in Mathematics) from 08.12.2015 to 06.02.2016.

The problem was formulated how to calculate so called the Mackenzie obstruction for existing of transitive Lie algebroid for given coupling between a finite dimensional Lie algebra and the tangent bundle of a smooth manifold.

It is proved that the Mackenzie obstruction for transitive Lie algebroids is trivial for the finite dimensional Lie algebra which is the direct sum of the center and the subalgebra without the center.

\end{abstract}

\section{Introduction.}

The preprint is prepared as description of results that were obtained during
joint scientific project No: 71NC /2015/VNCCCT on the VIASM (Vietnam Institute for Advanced Study in Mathematics) from 08.12.2015 to 06.02.2016.

Research results presented in this paper are a natural part of the overall program, which is carried out by a large group of researchers from different countries. This program essentially was initiated recently prematurely deceased Professor of Polytechnic University in Lodz (Poland), Jan Kubarski, who started  the joint work with one of the authors of this work to study the signatures of transitive Lie algebroids on compact manifolds (\cite{Kubarski-Mishchenko-2003},\cite{Kubarski-Mishchenko-2003-1},
\cite{Kubarski-Mishchenko-2004},\cite{Kubarski-Mishchenko-2009}).

The entire program can be described as the problem of homotopy classification of transitive Lie algebroids with fixed the base manifold and a fixed finite-dimensional Lie algebra as the associated Lie algebra for transitive Lie algebroids.

This task was actively engaged in collaboration consisting of the following mathematicians:
Vagif Ali-Muhtar Ogly Kasimov (Baku State University), Xiaoyu Li (Harbin Institute of Technology, China), A.S.Mishchenko (Moscow State Lomonosov University), Leanh Nguyen (National University of Civil Engineering, Vietnam), Jose Ribeiro Oliveira (Braga University, Portugal), Robert Wolak (Jagellonian University in Krakow, Poland)(\cite{Mishchenko-2010}-- \cite{Li-Mishchenko-2016}). Without the participation of the collaboration in the development of the project it would be impossible to solve a number of key issues and formulate remaining challenges.

\section{Formulation of the problem.}

 Let $\rg$ be a finite dimensional Lie algebra, let $\Aut(g)^{\delta}$  be the automorphism group of the Lie algebra $\rg$, with topology that on the quotient group  $\Out(g)$  the topology ia discrete. Consider a locally trivial bundle $L$ with typical fiber $\rg$ and structural group $\Aut(g)^{\delta}$. Fix on  $L$ correspondent local trivialization and associated with this trivialization the linear connection $\nabla$. The connection $\nabla$ defines the homomorphism od the bundles
$$
\nabla:TM\mapr{}\cD_{der}(L),
$$
 where $TM$ is the tangent bundle to manifold $M$,  $\cD_{der}(L)$ is the Lie algebroid of all covariant differentiations that satisfy Newton-Leibnitz relations with respect to fiberwise commutator operations (derivations) (so called bracket derivations).

Consider the exact sequence of Lie algebroids
\begin{equation}
0\mapr{}ZL\mapr{}L\mapr{ad}\cD_{der}(L)\mapr{\theta}\cD_{out}(L)\mapr{}0,
\end{equation}
the connection $\nabla$, associated with the trivialization in the bundle $L$ and induced by the homomorphism  $\theta$ connection
$$
\Xi:TM\mapr{}\cD_{out}(L), \quad \Xi=\theta\circ\nabla,
$$

$$
\xymatrix{
0\ar[r]&ZL\ar[r]&L\ar[r]^{ad}&\cD_{der}(L)\ar[r]^{\theta}&\cD_{out}(L)\ar[r]&0\\
&&&TM\ar[u]^{\nabla}\ar[ur]_{\Xi}&&
}
$$
he curvature tensir $R^\nabla$ satisfies the condition:
$$
\nabla\circ R^\nabla=R^\Xi=0.
$$

$$
\xymatrix{
0\ar[r]&ZL\ar[r]&L\ar[r]^{ad}&\Der(L)\ar[r]^{\theta}&\Out(L)\ar[r]&0\\
&&&\Lambda^{2}TM\ar[u]^{R^{\nabla}}\ar[ur]_{R^{\Xi}=0}\ar[lu]^{\Omega}&&
}
$$
Hence,  due to exactness of the sequence (1) we have:
$$
R^\nabla=ad\circ\Omega
$$
foe some form  $\Omega$.

Then the differential $d^{\nabla}\Omega$ of the form $\Omega$ can be included in the commutative diagram

$$
\xymatrix{
0\ar[r]&ZL\ar[r]^{i}&L\ar[r]^{ad}&\Der(L)\ar[r]^{\theta}&\Out(L)\ar[r]&0\\
&&\Lambda^{3}TM\ar[ru]_{d^{\nabla^{der}}R^{\nabla}=0}\ar[u]^{d^{\nabla}\Omega}\ar[lu]^{d^{\nabla}\Omega}&&
}
$$
hence values in $ZL$:

$$
\xymatrix{
\Lambda^{3}TM\ar[r]^{d^{\nabla}\Omega}& ZL.
}
$$

Since the connection $\nabla$ make the subbundle $ZL$ to be invariant and induces on the subbundle $ZL$ the connection $\nabla^{Z}$ with trivial curvature tensor,
$$
R^{\nabla^{Z}}\equiv 0,
$$
such that
$$
d^{\nabla^{Z}}\left(d^{\nabla}\Omega\right)\equiv 0,
$$
the form $d^{\nabla}\Omega$ is closed form in the space
$\Omega^{3}\left(M;ZL;\nabla^{Z}\right)$. The cohomology class
$H^{3}\left(M;ZL;\nabla^{Z}\right)$, induced by the form $d^{\nabla}\Omega$,
is denoted by
$$
\cO bs(\nabla, \Omega)\in H^{3}\left(M;ZL;\nabla^{Z}\right)
$$
and is called the Mackenzie obstruction for the coupling $\Xi$.

Due to the theorem 7.2.12 from \cite{Mackenzie-05} (p. 277) The cohomology class
$\cO bs(\nabla, \Omega)\in H^{3}\left(M;ZL;\nabla^{Z}\right)$ depends only on coupling
$\Xi$ and does not depend of the choice of the connection $\nabla$ nd os the form $\Omega$, that is
$$
\cO bs(\nabla, \Omega)=\cO bs(\Xi).
$$

This obstruction is trivial then and only then when the pair $(\nabla, \Omega)$ can be chosen such that  $d^{\nabla}\Omega\equiv 0$, that is the form $\Omega$ forma a structure of a transitive Lie algebroid for given coupling $\Xi$.

So, for classification of transitive Lie algebroids with fixed adjoint finite dimensional Lie algebra $\rg$ it is necessary to check for which coupling $\Xi$ the cohomology class $\cO bs(\Xi)\in H^{3}\left(M;ZL;\nabla^{Z}\right)$ is trivial. Therefore a natural problem arises about calculation and description of the cohomology class
$\cO bs(\Xi)\in H^{3}\left(M;ZL;\nabla^{Z}\right)$.

First, note that each coupling $\Xi$ is in one-to-one correspondence
with the class of locally trivial bundle structures $L$ with structure
group $\Aut (g)^{\delta}$, or
that is the same, with the homotopy classes of maps of the base $M$ to the classifying space
$B(\Aut (g)^{\delta})$,
$$\begin{array}{l}
\{\Xi:TM\mapr{}\cD_{out}(L)\}\Leftrightarrow [M,B(\Aut(g)^{\delta})],\\\\
\Xi=\Xi[f] \hbox{ for }
f:M\mapr{}B(\Aut(g)^{\delta}).
\end{array}
$$

This correspondence naturally generates the assumption that an obstruction to coupling of
$\cO bs(\Xi)\in H^{3}\left(M; ZL; \nabla^{Z}\right)$ is a characteristic class
for the bundle $L$ with the structural group $\Aut(g)^{\delta}$ with values in the cohomology
$H^{3}\left(M; ZL; \nabla^{Z}\right)$. If so, then the Mackenzie obstruction for the
coupling of $\cO bs(\Xi)$
can be calculated as the inverse image under a continuous map $ f:M\mapr{}B(\Aut(g)^{\delta})$
of a cohomology class
$$
\cO bs^{\infty}\in H^{3}\left(B(\Aut(g)^{\delta});ZL_{B(\Aut(g)^{\delta})};\nabla^{Z}\right),
$$
$$
\cO bs(\Xi)=\cO bs(\Xi[f])=f^{*}(\cO bs^{\infty}).
$$

To prove that the Mackenzie obstruction for the coupling of
$$\cO bs(\Xi)\in H^{3}\left(M;ZL;\nabla^{Z}\right)$$
is the characteristic class
for the bundle $L$ with structure group $\Aut(g)^{\delta}$ with values in the cohomology
$H^{3}\left(M;ZL;\nabla^{Z}\right)$ requires the following assertion:

\begin{theorem}{}
Consider two manifolds $M_{1}$ and $M_{2}$ and a smooth map $f:M_{1}\mapr{}M_{2}$.
Let $L_{2}$ - locally trivial bundle with fiber $\rg$ and structure group $\Aut(g)^{\delta}$ on the base $M_{2}$. Let $L_{1}=f^{*}(L_{2})$ be the  pullback bundle
$L_{2}$. Let $\Xi_{1}$ and $\ Xi_{2}$  be couplings for the corresponding bundle $L_{1}$
and $L_{2}$.
Then
$$
\cO bs(\Xi_{1}) = f^{*}(\cO bs(\Xi_{2}))
$$
where $f^{*}$ is homomorphism in cohomologies:
$$
f^{*}: H^{3}(M_{2};ZL_{2};\nabla^{Z})\mapr{}H^{3}(M_{1};ZL_{1};\nabla^{Z}).
$$
\end{theorem}
This theorem was firstly proved by Li Xiaoyu (unpublished).

Based on the above theorem the stated problem of calculating and describing the Mackenzie obstructions to coupling can be reduced to the description of the cohomology class
$$
\cO bs^{\infty}\in H^{3}\left(B(\Aut(g)^{\delta}); ZL_{B(\Aut(g)^{\delta})}; \nabla^{Z}\right),
$$
in particular, the finding where the class is equal to zero. In the book by the Mackenzie this issue was not discussed. In other words, it was not discussed if the Mackenzie obstructions for coupling is not nontrivial.
This work is devoted to just the question whether or not the Mackenzie obstruction for coupling is trivial.

\section {The simplest examples.} The question of calculating of the Mackenzie obstructions for coupling can be considered
first for the simplest examples of finite-dimensional Lie algebras. Two extreme cases are the most simple.

The first extreme case where the Lie algebra $\rg$ has no center, $Z\rg = 0$.
In this case, the Mackenzie obstruction for coupling is trivial.

The second extreme case where the Lie algebra coincides with its center $\rg = Z\rg $, i.e. the Lie algebra is commutative. In this case, the group $\Aut(g)^{\delta}$ is isomorphic to the group of all matrices $GL(\rg)$ with discrete topology. In this case, the coupling $\Xi$ coincides with a flat connection $\nabla$ in a flat bundle $L$, i.e. $R^{\nabla} \equiv 0 $. This means that the form $\Omega$ can be chosen trivial, i.e. $d^{\nabla}\Omega = 0 $. So the Mackenzie obstruction for coupling of $\cO bs(\Xi)$
equals to zero.

Third , the general case where $Z\rg\subset\g$, $Z\rg\neq\rg$.

As a general hypothesis, we conjecture that the Mackenzie obstruction to the coupling is always zero. Even if this hypothesis is not correct,  further consideration will help to define more precisely the conditions under which the coupling of an the Mackenzie obstruction to be nontrivial.

\section {The reduction to the quotient algebra  by its center}

Consider a finite-dimensional Lie algebra $\rg$, in which the center $Z\rg$ does not coincide with the algebra
$\rg$. Let $\ rg_{0} = \rg/Z\rg $, i.e. is an exact sequence
$$
0\mapr{}Z\rg\mapr{}\rg\mapr{}\rg_{0}\mapr{}0.
$$

It is easy to check that every automorphism $\varphi\in\Aut(\rg)$ leaves the center $ Z\rg$ invariant, i.e. induces an automorphism of the quotient algebra $\rg_{0}$. Consequently, it turns out
natural homomorphism
$$
\Aut(\rg)\mapr{\alpha}\Aut(\rg_{0}),
$$
which takes the subgroup $\Int(\rg)$ to subgroup $\Int(\rg_{0})$.
Hence the homomorphism $\alpha$ is  well-defined homomorphism
$$
\Aut(\rg)^{\delta}\mapr{\alpha}\Aut(\rg_{0})^{\delta}.
$$
This observation allows to state the theorem:

\begin{theorem}
For the bundle LAB with the Lie algebra
$\rg_{0}$ as the fiber and the structural group $\Aut(\rg_{0})^{\delta}$ that can be reduced to the group $\Aut(\rg)^{\delta}$ the Mackenzie obstruction of the LAB is trivial.
\end{theorem}

\section {The case of direct sum of the quotient algebra and the center}

There is another case how to analyze the Mackenzie obstruction for given coupling
between the Lie algebra bundle $L$ (LAB) and the tangent bundle $TM$
when $g=Zg\oplus g_{0}$, $Zg_{0}=0$.

If

$$
\left(
\begin{array}{c}
z_{1}\\u_{1}
\end{array}\right),
\left(
\begin{array}{c}
z_{2}\\u_{2}
\end{array}\right)\in g
$$
then

$$
\left[
\left(
\begin{array}{c}
z_{1}\\u_{1}
\end{array}\right),
\left(
\begin{array}{c}
z_{2}\\u_{2}
\end{array}\right)
\right]=
\left(
\begin{array}{c}
0\\
\left[u_{1},u_{2}\right]
\end{array}\right)
\in g
$$
Consider the automorphism group $\Aut(g)$ as a family of matrices
$$\varphi=
\left(\begin{array}{cc}
\varphi^{1}_{1} & \varphi^{1}_{2}\\
\varphi^{2}_{1} & \varphi^{2}_{2}
\end{array}\right):
Zg\oplus g_{0}\mapr{}Zg\oplus g_{0}.
$$
The automorphism $\varphi$ should satisfy the condition:
$$
\varphi\left(
\left[
\left(
\begin{array}{c}
z_{1}\\u_{1}
\end{array}\right),
\left(
\begin{array}{c}
z_{2}\\u_{2}
\end{array}\right)
\right]
\right)=
\left[
\varphi\left(
\begin{array}{c}
z_{1}\\u_{1}
\end{array}\right),
\varphi\left(
\begin{array}{c}
z_{2}\\u_{2}
\end{array}\right)
\right],
$$

that is

$$
\left(\begin{array}{cc}
\varphi^{1}_{1} & \varphi^{1}_{2}\\
\varphi^{2}_{1} & \varphi^{2}_{2}
\end{array}\right)
\left(
\left[
\left(
\begin{array}{c}
z_{1}\\u_{1}
\end{array}\right),
\left(
\begin{array}{c}
z_{2}\\u_{2}
\end{array}\right)
\right]
\right)=
\left[
\left(\begin{array}{cc}
\varphi^{1}_{1} & \varphi^{1}_{2}\\
\varphi^{2}_{1} & \varphi^{2}_{2}
\end{array}\right)
\left(
\begin{array}{c}
z_{1}\\u_{1}
\end{array}\right),
\left(\begin{array}{cc}
\varphi^{1}_{1} & \varphi^{1}_{2}\\
\varphi^{2}_{1} & \varphi^{2}_{2}
\end{array}\right)
\left(
\begin{array}{c}
z_{2}\\u_{2}
\end{array}\right)
\right],
$$

or

$$
\left(\begin{array}{cc}
\varphi^{1}_{1} & \varphi^{1}_{2}\\
\varphi^{2}_{1} & \varphi^{2}_{2}
\end{array}\right)
\left(
\begin{array}{c}
0\\
\left[u_{1},u_{2}\right]
\end{array}
\right)=
\left[
\left(\begin{array}{cc}
\varphi^{1}_{1} & \varphi^{1}_{2}\\
\varphi^{2}_{1} & \varphi^{2}_{2}
\end{array}\right)
\left(
\begin{array}{c}
z_{1}\\u_{1}
\end{array}\right),
\left(\begin{array}{cc}
\varphi^{1}_{1} & \varphi^{1}_{2}\\
\varphi^{2}_{1} & \varphi^{2}_{2}
\end{array}\right)
\left(
\begin{array}{c}
z_{2}\\u_{2}
\end{array}\right)
\right],
$$

or

$$
\left(\begin{array}{c}
\varphi^{1}_{2}(\left[u_{1},u_{2}\right])\\
\varphi^{2}_{2}(\left[u_{1},u_{2}\right])
\end{array}\right)
=
\left[
\left(\begin{array}{c}
\varphi^{1}_{1}(z_{1}) + \varphi^{1}_{2}(u_{1})\\
\varphi^{2}_{1}(z_{1}) + \varphi^{2}_{2}(u_{1})
\end{array}\right)
,
\left(\begin{array}{cc}
\varphi^{1}_{1}(z_{2}) + \varphi^{1}_{2}(u_{2})\\
\varphi^{2}_{1}(z_{2}) + \varphi^{2}_{2}(u_{2})
\end{array}\right)
\right],
$$

or

$$
\left(\begin{array}{c}
\varphi^{1}_{2}(\left[u_{1},u_{2}\right])\\
\varphi^{2}_{2}(\left[u_{1},u_{2}\right])
\end{array}\right)
=
\left(\begin{array}{c}
0\\
\left[
\varphi^{2}_{1}(z_{1}) + \varphi^{2}_{2}(u_{1}),
\varphi^{2}_{1}(z_{2}) + \varphi^{2}_{2}(u_{2})
\right]
\end{array}\right)
,
$$

Hence if $z_{1}=z_{2}=0$, then

$$
\left(\begin{array}{c}
\varphi^{1}_{2}(\left[u_{1},u_{2}\right])\\
\varphi^{2}_{2}(\left[u_{1},u_{2}\right])
\end{array}\right)
=
\left(\begin{array}{c}
0\\
\left[
\varphi^{2}_{2}(u_{1}),
\varphi^{2}_{2}(u_{2})
\right]
\end{array}\right)
,
$$

that is $\varphi^{2}_{2}\in \Aut(g_{0})$.

If $u_{1}=0$ then
$$
\left[
\varphi^{2}_{1}(z_{1}),
\varphi^{2}_{1}(z_{2}) + \varphi^{2}_{2}(u_{2})
\right]=0,
$$
that is
$$
\varphi^{2}_{1}(z_{1})=0
$$
or

$$
\varphi^{2}_{1}=0.
$$

Thus for any $\varphi\in\Aut(g)$ one has
$$
\varphi=
\left(\begin{array}{cc}
\varphi^{1}_{1} & \varphi^{1}_{2}\\
0 & \varphi^{2}_{2}
\end{array}\right):
Zg\oplus g_{0}\mapr{}Zg\oplus g_{0},
$$

$$
\varphi^{2}_{2}\in \Aut(g_{0}),
$$

$$
\varphi^{1}_{2}|_{[g_{0},g_{0}]}=0.
$$

Hence
$$
\Aut(g)=\left(\begin{array}{cc}
\GL(Zg)& \Hom(g_{0}/[g_{0},g_{0}], Zg_{0})\\
0&\Aut(g_{0})
\end{array}\right),
$$
where all items are independent each from others

Let us describe $\Inn(g)\subset\Aut(g)$. Each element in $\varphi\in\Inn(g)$ has the representation
$$
\varphi=\exp\ad(\sigma), \sigma\in g,
$$
where $\ad(\sigma)(\sigma')=[\sigma,\sigma']$.

Let $\sigma=\left(\begin{array}{c}
z \\ u
\end{array}\right)$.
Then
$$
\ad(\sigma)=\left(
\begin{array}{cc}
0&0\\
0&\ad(u)
\end{array}
\right).
$$
Hence
$$
\exp\ad(\sigma)=
\left(
\begin{array}{cc}
1&0\\
0&\exp\ad(u)
\end{array}
\right).
$$
So the group $\Inn(g)$ consists of the matrices
$$
\Inn(g)=\left(\begin{array}{cc}
1&0\\
0&\Inn(g_{0})
\end{array}\right)
$$
Hence the quotient group $\Aut(g)/\Inn(g)$ is
$$
\Aut(g)/\Inn(g)=
\left(\begin{array}{cc}
\GL(Zg)& \Hom(g_{0}/[g_{0},g_{0}], Zg)\\
0&\Aut(g_{0})/\Inn(g_{0})
\end{array}\right).
$$
From the  result one can show that in this case the Mackenzie obstruction is trivial.

\vskip 1cm
{\Large\bf The publications of the collaboration:}


\begin{thebibliography}{10}

\bibitem{Mackenzie-05}
K.C.H. Mackenzie.
\newblock {\em General Theory of Lie Groupoids and Lie Algebroids}.
\newblock Cambridge University Press, 2005.
\end{thebibliography}

\begin{thebibliography}{aaa}

\bibitem{Kubarski-Mishchenko-2003}
J.Kubarski, A.S. Mishchenko,
{\em Transitive Lie algebroids: spectral sequences and signature,}
5th International Conference Geometry and Topology of Manifolds, Krynica-Zdruj, Poland, 2003

\bibitem{Kubarski-Mishchenko-2003-1}
J.Kubarski, A.S. Mishchenko,
{\em Lie algebroids: Spectral sequences and signature,}
Doklady Mathematics, Maik Nauka/Interperiodica Publishing (Russian Federation), t. 68, No.2, p. 188-190

\bibitem{Kubarski-Mishchenko-2004}
J.Kubarski, A.S. Mishchenko,
{\em Nondegenerate cohomology pairing for transitive Lie algebroids, characterization,}
Central European Journal of Mathematics, Springer Verlag (Germany),2, No.5, p. 663-707

\bibitem{Kubarski-Mishchenko-2009}
J.Kubarski, A.S. Mishchenko,
{\em Algebraic aspects of the Hirzebruch signature operator and applications to transitive Lie algebroids,}
Russian Journal of Mathematical Physics, Maik Nauka/Interperiodica Publishing (Russian Federation), 16,3, p. 413-428, 2009

\bibitem{Mishchenko-2010}
A.S. Mishchenko,
{\em Transitive Lie algebroids - categorical point of
view,}
arXiv:submit/0064778 [math.AT] 24 Jun 2010

\bibitem{Mishchenko-Ribeiro-2011}
A.S. Mishchenko, J.Ribeiro
{\em Generalization of the Sullivan construction for
Transitive Lie Algebroids,}
arXiv:submit/0370105 [math.AT] 29 Nov 2011



\bibitem{Li-Mishchenko-2012}
Xiaoyu Li,
A.S. Mishchenko
{\em Comparison of categorical characteristic classes of transitive Lie algebroid with Chern-Weil homomorphism,}
Electronic archive: arXiv:1208.6564, 2012

\bibitem{Li-Mishchenko-2012-1}
Xiaoyu Li,
A.S. Mishchenko
{\em Some geometrical problems of construction of the classifying space of the transitive Lie algebroid,}
The second International Conference "K-Theory, C*-algebras and Topology of Manifolds II", Tianjin, China, 2012


\bibitem{Li-Mishchenko-2013}
Xiaoyu Li,
A.S. Mishchenko
{\em The Existence of Coupling in the Category of
Transitive Lie Algebroid,}
arXiv:submit/0744505 [math.AT] 23 Jun 2013

\bibitem{Li-Mishchenko-2013-1}
Xiaoyu Li,
A.S. Mishchenko
{\em Description of coupling in the category of
transitive Lie algebroids,}
arXiv:1310.5824v1 [math.AT] 22 Oct 2013

\bibitem{Li-Mishchenko-Gasimov-2014}
Xiaoyu Li,
A.S. Mishchenko, V.Gasimov
{\em Mackenzie obstruction for the existence of a transitive Lie algebroid,}
 Russian Journal of Mathematical Physics, Maik Nauka/Interperiodica Publishing (Russian Federation), t. 21, No.4, p. 544-548, 2014

 \bibitem{Li-Mishchenko-2015}
Xiaoyu Li,
A.S. Mishchenko
{\em Classification of Couplings for Transitive Lie Algebroids,}
 Doklady Mathematics, Maik Nauka/Interperiodica Publishing (Russian Federation), t. 91, No.1, p. 84-86, 2015

\bibitem{Li-Mishchenko-2016}
Xiaoyu Li,
A.S. Mishchenko
{\em The existence and classification of couplings between Lie algebrabundles and tangent bundles,}
Topology and its Applications, Elsevier BV (Netherlands), t. 200, p. 1-18, 2016

\end{thebibliography}
\end{document}